\documentclass{amsart}
\usepackage{mathrsfs}
\usepackage{amsfonts}
\usepackage{amssymb}
\usepackage[papersize={6.9in,9.8in},textwidth=14cm,textheight=21.7cm,centering]{geometry}

\usepackage[colorlinks=true,citecolor=blue,linkcolor=blue]{hyperref}

\usepackage{hyperref}

\allowdisplaybreaks

\newtheorem*{theorem*}{Theorem}

\newtheorem{theorem}{Theorem}

\newtheorem{lemma}{Lemma}

\newtheorem*{proposition*}{Proposition}

\newtheorem*{corollary*}{Corollary}

\theoremstyle{remark}

\newcommand{\ud}{\mathrm{d}}

\begin{document}

\title[Moments of Kloosterman sums]{A note on the moments of Kloosterman sums}

\author{Ping Xi}

\author{Yuan Yi}

\address{School of Science, Xi'an Jiaotong University, Xi'an 710049, P. R. China}

\email{xprime@163.com}



\begin{abstract}In this note, we deduce an asymptotic formula for even power moments of Kloosterman sums based on the important work of N. M. Katz on Kloosterman sheaves. In a similar manner, we can also obtain an upper bound for odd power moments. Moreover, we shall give an asymptotic formula for odd power moments of absolute Kloosterman sums. Consequently, we find that there are infinitely many $a\bmod p$ such that $S(a,1;p)\gtrless0$ as $p\rightarrow+\infty.$ \end{abstract}

\maketitle

\section{Introduction}

Let $p$ be an odd prime. For any integer $a$ and $b$, the classical Kloosterman sum is defined by
\[S(a,b;p)=\sideset{}{^*}\sum_{x\bmod p}e\left(\frac{ax+b\overline{x}}{p}\right),\]where the summation is over a reduced residue system mod $p$. Such an exponential sum first appeared in a paper of H. Poincar\'{e} \cite{Po} on modular functions.
As tackling the problem on the representation of numbers in positive definite diagonal quadratic forms, H. D. Kloosterman \cite{Kl} re-introduced and first seriously studied the exponential sum $S(a,b;p)$, which was later named after him. He proved for any odd prime $p$ that
\begin{align}\label{eq:1}|S(a,b;p)|\leqslant3^{1/4}p^{3/4},\ \ \ (ab,p)=1.\end{align} The best result up to now is due to A. Weil \cite{We}
\begin{align}\label{eq:2}|S(a,b;p)|\leqslant2p^{1/2},\ \ \ (ab,p)=1,\end{align} the proof of which originally requires deep results of algebraic geometry. At present, due to the important work of S. A. Stepanov, W. Schmidt and E. Bombieri, Weil's bound can be established by elementary means. See \cite{Sch} for details.

Kloosterman sums play quite an important role in modern analytic number theory, as well as the theory of automorphic forms (see \cite{Iw} for instance). One of the problems is to consider the moments of Kloosterman sums, namely
\[V_k(p)=\sideset{}{^*}\sum_{a\bmod p}S^k(a,1;p).\] The original motivation to consider the moments just follows H. D. Kloosterman, who wished to gain a nontrivial bound for the individual sum. This leads him to consider the problem in a global sense first. For the cases $k\leqslant4$, we can give exact identities by elementary methods. To be precise, we have
\begin{align*}&V_1(p)=1,\\
&V_2(p)=p^2-p-1,\\
&V_3(p)=\bigg(\frac{p}{3}\bigg)p^2+2p+1,\\
&V_4(p)=2p^3-3p^2-3p-1,\end{align*}
where $(\frac{\cdot}{p})$ is the Legendre symbol mod $p$. The last one was first obtained by H. D. Kloosterman, which enables him to deduce the first nontrivial bound (\ref{eq:1}).

For the case of $k=5$, it follows from the work in \cite{Li,PTV} that
\[V_5(p)=\bigg(\frac{p}{3}\bigg)4p^3+(a_p+5)p^2+4p+1,\ \ p>5,\]where
where $a_p$ is the integer with $|a_p|<2p$ defined for $p>5$ by
\[a_p=\begin{cases}
2p-12u^2,\ &\text{ if }p=3u^2+5v^2,\\
4x^2-2p,\ &\text{ if }p=x^2+15y^2,\\
0,\ &\text{ if }p\equiv7,11,13 \text{ or }14\pmod{15}.\end{cases}\]

For $k=6$, H. Sali\'{e} \cite{Sa} and H. Davenport \cite{Da} independently proved in an elementary manner that
\[V_6(p)\ll p^4,\]which leads to the result that the exponent in (\ref{eq:1}) can be reduced from $3/4$ to $2/3$ up to the constant factor. On the other hand, it follows from the work in \cite{HSvv} that
\[V_6(p)=5p^4-10p^3-(b_p+9)p^2-5p-1,\ \ p>7,\]
where $b_p$ is the integer with $|b_p|<2p^{3/2}$  defined to be the $p$-th Fourier coefficient in the Fourier expansion of the newform of weight 4,
level 6 given by $(\eta(z)\eta(2z)\eta(3z)\eta(6z))^2$, here $\eta$ is the Dedekind eta function.

In a recent paper, R. J. Evans \cite{Ev} studied the case of $k=7$ and conjecturely obtained an exact identity in terms of Hecke eigenvalues for a weight 3 newform on $\Gamma_0(525)$ with quartic nebentypus of conductor 105.

However, for the cases $k\geqslant8$, it seems that there are not too many works focusing on the identities of moments of Kloosterman sums. It should be remarked that N. M. Katz \cite{Ka} has studied the higher moments from a modern point of view. From his result, one can deduce certain asymptotic formulae for any even power moments of $S(a,1;p)$. However, we can just obtain upper bounds for the case of odd powers, and it seems that to get an asymptotic formula valid for a general odd integer $k$ is out of reach at present.

\begin{theorem}\label{thm:1} For any given natural number $k$ and each large prime number $p,$ we have
\begin{align*}V_{2k}(p)=\frac{1}{k+1}\binom{2k}{k}p^{k+1}+O(p^{k+1/2}),\end{align*}
\begin{align*}V_{2k+1}(p)\ll p^{k+1},\end{align*}where $\binom{m}{n}$ is the binomial coefficient defined as $m!/n!(m-n)!,$ and the implied constants depend only on $k.$
\end{theorem}

In fact, if we apply Weil's bound (\ref{eq:2}) to each individual Kloosterman sum and then summing over $a$, we can deduce the rough estimate
\[|V_{2k+1}(p)|\leqslant\sideset{}{^*}\sum_{a\bmod p}|S(a,1;p)|^{2k+1}\ll p^{k+3/2}.\]Hence one can expect that there must be significant cancellations among these individual Kloosterman sums in $V_{2k+1}(p)$. In order to illustrate this phenomenon, we would like to investigate the following moments of \emph{absolute} Kloosterman sums:
\[\widetilde{V}_{k}(p)=\sideset{}{^*}\sum_{a\bmod p}|S(a,1;p)|^{k}.\] In fact, we can establish an asymptotic formula for $\widetilde{V}_{k}(p)$ as $p$ tends to infinity.
\begin{theorem}\label{thm:2} For any given natural number $k$ and each large prime number $p,$ we have
\begin{align*}\widetilde{V}_{2k+1}(p)=\frac{4^{k+1}}{\pi(2k+1)(2k+3)}\binom{2k}{k}^{-1}p^{k+3/2}+O_k(p^{k+1}),\end{align*}
where the implied constant depends only on $k.$
\end{theorem}

Subsequently, we can deduce from Theorems \ref{thm:1} and \ref{thm:2} that
\begin{corollary*}\label{coro:1} For sufficiently large prime number $p,$ we have
\begin{align*}\sideset{}{^*}\sum_{\substack{a\bmod p\\S(a,1;p)>0}}1\geqslant\frac{4}{9\pi^2}p+O(p^{1/2}),\end{align*}and
\begin{align*}\sideset{}{^*}\sum_{\substack{a\bmod p\\S(a,1;p)<0}}1\geqslant\frac{4}{9\pi^2}p+O(p^{1/2}).\end{align*}

\end{corollary*}

We shall give the proof of the theorems along the line of the following sections. First, we shall review some basic properties of Chebyshev polynomials in Section \ref{sec:2}, and then we shall deduce the theorems from Katz's result on Kloosterman sheaves, which will be completed in Section \ref{sec:3}. The proof of the corollary will be completed in Section \ref{sec:4}.

\bigskip

\section{Chebyshev polynomials}\label{sec:2}

The Chebyshev polynomials $U_k(x)(k\geqslant0)$ are defined recursively by
\[U_0(x)=1,\ \ U_1(x)=2x,\]
\[U_{k+1}(x)=2xU_k(x)-U_{k-1}(x).\]
Chebyshev polynomials enjoy the following identity
\[\frac{2}{\pi}\int_{-1}^1\sqrt{1-x^2}U_m(x)U_n(x)\ud x=\begin{cases}
1,\ &m=n,\\
0,\ &m\neq n.\end{cases}\]
This shows that Chebyshev polynomials are orthogonal with respect to the weight
\[\frac{2}{\pi}\sqrt{1-x^2}\]
on the interval $[-1,1]$. Hence one may expect that every continuous function defined in $[-1,1]$ can be represented by a certain linear combination of such Chebyshev polynomials. In particular, we have
\begin{lemma}[\cite{Ri}, P.54] For each $k\in\mathbb{N}$ and $x\in[-1,1]$, we have
\begin{align}\label{eq:3}x^{2k}=\frac{(2k)!}{4^k}\sum_{0\leqslant\ell\leqslant k}\frac{2\ell+1}{(k-\ell)!(k+\ell+1)!}U_{2\ell}(x),\end{align}
and
\begin{align}\label{eq:4}x^{2k+1}=\frac{(2k+1)!}{4^{k+1}}\sum_{0\leqslant\ell\leqslant k}\frac{\ell+1}{(k-\ell)!(k+\ell+2)!}U_{2\ell+1}(x).\end{align}
\end{lemma}

Moreover, we can express $|x|^{2k+1}$ as a linear combination of certain Chebyshev polynomials as follows:
\begin{align*}|x|^{2k+1}=\sum_{\ell\geqslant0}a_\ell U_{\ell}(x),\end{align*}where the coefficients $a_j$'s could be computed in the following manner:
\begin{align*}a_\ell=\frac{2}{\pi}\int_{-1}^1\sqrt{1-x^2}U_\ell(x)|x|^{2k+1}\ud x.\end{align*}

First we have
\begin{align*}a_\ell&=\frac{2}{\pi}\int_0^1\sqrt{1-x^2}x^{2k+1}(U_\ell(x)+U_\ell(-x))\ud x.\end{align*}
Putting $x=\cos\theta$, we obtain that
\begin{align*}a_\ell&=\frac{2}{\pi}\int_{\pi/2}^{0}\sin\theta(\cos\theta)^{2k+1}(U_\ell(\cos\theta)+U_\ell(-\cos\theta))\ud\cos\theta\\
&=\frac{2}{\pi}\int_{0}^{\pi/2}\sin^2\theta(\cos\theta)^{2k+1}(U_\ell(\cos\theta)+U_\ell(-\cos\theta))\ud\theta.\end{align*}
Since
\[U_\ell(\cos\theta)=\frac{\sin((\ell+1)\theta)}{\sin\theta},\ \ U_\ell(-\cos\theta)=(-1)^{\ell}\frac{\sin((\ell+1)\theta)}{\sin\theta},\]
it follows that
\begin{align*}a_\ell&=\frac{2}{\pi}(1+(-1)^{\ell})\int_{0}^{\pi/2}\sin\theta(\cos\theta)^{2k+1}\sin((\ell+1)\theta)\ud\theta.\end{align*}
Thus $a_\ell$=0 if $\ell$ is odd, and if $\ell$ is even, we have
\begin{align}\label{eq:5}a_\ell&=\frac{4}{\pi}\int_{0}^{\pi/2}\sin\theta(\cos\theta)^{2k+1}\sin((\ell+1)\theta)\ud\theta.\end{align}

For $\ell=0$, we have
\begin{align*}a_0&=\frac{4}{\pi}\int_{0}^{\pi/2}\sin^2\theta(\cos\theta)^{2k+1}\ud\theta=\frac{4^{k+1}}{\pi(2k+1)(2k+3)}\binom{2k}{k}^{-1}.\end{align*}
For $\ell\geqslant2k+2$, we have
\begin{align*}a_\ell&=\frac{2}{\pi}\int_{0}^{\pi/2}(\cos(\ell\theta)-\cos((\ell+2)\theta))(\cos\theta)^{2k+1}\ud\theta,\end{align*}
then it follows that (by \cite{RG}, 2.538)
\[a_\ell=\frac{2}{\pi}\left(\dfrac{(-1)^{\ell/2+k+1}(2k+1)!}{(2k+1+\ell)(\ell-2k-1)!}-\dfrac{(-1)^{\ell/2+k}(2k+1)!}{(2k+3+\ell)(\ell-2k+1)!}\right).\]
Moreover, we have
\[a_\ell\ll\frac{1}{(\ell-2k)!}\]for any $\ell\geqslant2k+2$, where the implied constant depends only on $k$.

Hence we can conclude that
\begin{lemma}\label{lm:2}For each $k\in\mathbb{N}$ and $x\in[-1,1]$, we have
\begin{align}\label{eq:6}|x|^{2k+1}=\frac{4^{k+1}}{\pi(2k+1)(2k+3)}\binom{2k}{k}^{-1}+\sum_{\ell\geqslant1}c_{\ell,k} U_{2\ell}(x),\end{align}where $c_{\ell,k}=a_{2\ell}$ could be computed exactly by $(\ref{eq:5})$ for $1\leqslant\ell\leqslant k$ and is bounded as
\[c_{\ell,k}\ll_k\frac{1}{(2\ell-2k)!}\]for $\ell\geqslant k+1.$
\end{lemma}

\bigskip

\section{Derivation from Katz's estimate}\label{sec:3}

In view of Weil's bound (\ref{eq:2}), we can write
\begin{align}\label{eq:7}\frac{S(a,1;p)}{\sqrt{p}}=2\cos\theta_p(a),\end{align}where $0\leqslant\theta_p(a)\leqslant\pi$ is called the Kloosterman sum angle.

The main result of N. M. Katz \cite{Ka} can be stated equivalently as follows.
\begin{proposition*}[see \cite{Iw}, Theorem 4.6]For any positive integer $k$ and odd prime $p$, we have
\begin{align}\label{eq:8}\bigg|~\sideset{}{^*}\sum_{a\bmod p}U_k(\cos\theta_p(a))\bigg|\leqslant\frac{1}{2}(k+1)p^{1/2}.\end{align}
\end{proposition*}

Following the notation in (\ref{eq:7}), we have
\begin{align*}V_{2k}(p)&=(4p)^{k}\sideset{}{^*}\sum_{a\bmod p}(\cos\theta_p(a))^{2k},\end{align*}from which and (\ref{eq:3}), we can deduce that
\begin{align*}V_{2k}(p)&=p^{k}(2k)!\sum_{0\leqslant\ell\leqslant k}\frac{2\ell+1}{(k-\ell)!(k+\ell+1)!}\sideset{}{^*}\sum_{a\bmod p}U_{2\ell}(\cos\theta_p(a)).\end{align*}

The term $\ell=0$ gives the contribution
\begin{align}\label{eq:9}p^{k}\frac{(2k)!}{k!(k+1)!}\sideset{}{^*}\sum_{a\bmod p}U_{0}(\cos\theta_p(a))=\frac{1}{k+1}\binom{2k}{k}p^{k+1}+O(p^k),\end{align}
where the $O$-constant depends only on $k$. For the terms with $1\leqslant\ell\leqslant k$, we can apply Katz's estimate (\ref{eq:8}) to each one, getting
\begin{align}\label{eq:10}p^{k}(2k)!\sum_{1\leqslant\ell\leqslant k}\frac{2\ell+1}{(k-\ell)!(k+\ell+1)!}\sideset{}{^*}\sum_{a\bmod p}U_{2\ell}(\cos\theta_p(a))\ll_k p^{k+1/2}.\end{align}

Combining (\ref{eq:9}) and (\ref{eq:10}), we finally arrive at
\begin{align*}V_{2k}(p)=\frac{1}{k+1}\binom{2k}{k}p^{k+1}+O(p^{k+1/2}),\end{align*}
where the $O$-constant depends only on $k$.

Following a similar argument, we can write
\begin{align*}V_{2k+1}(p)&=\frac{1}{2}(2k+1)!p^{k+1/2}\sum_{0\leqslant\ell\leqslant k}\frac{\ell+1}{(k-\ell)!(k+\ell+2)!}\sideset{}{^*}\sum_{a\bmod p}U_{2\ell+1}(\cos\theta_p(a)).\end{align*}Applying Katz's estimate (\ref{eq:8}), we can find that
\begin{align*}V_{2k+1}(p)\ll p^{k+1}.\end{align*}

This completes the proof of Theorem \ref{thm:1}. And Theorem \ref{thm:2} follows in a similar way, since for any fixed $k$, the coefficient $c_{\ell,k}$ in Lemma \ref{lm:2} decays rapidly as $\ell$ tends to infinity.

\bigskip

\section{Proof of the corollary}\label{sec:4}

Observing that
\[\sideset{}{^*}\sum_{\substack{a\bmod p\\S(a,1;p)>0}}S(a,1;p)=\frac{1}{2}(V_{1}(p)+\widetilde{V}_{1}(p)),\] thus we can conclude from Theorems \ref{thm:1} and \ref{thm:2} that
\begin{align}\label{eq:11}\sideset{}{^*}\sum_{\substack{a\bmod p\\S(a,1;p)>0}}S(a,1;p)=\frac{1}{3\pi}p^{3/2}+O(p).\end{align}

On the other hand, from Cauchy inequality we have
\begin{align*}\bigg(\sideset{}{^*}\sum_{\substack{a\bmod p\\S(a,1;p)>0}}S(a,1;p)\bigg)^2\leqslant V_{2}(p)\sideset{}{^*}\sum_{\substack{a\bmod p\\S(a,1;p)>0}}1,\end{align*}
it follows from (\ref{eq:11}) and Theorem \ref{thm:1} that
\begin{align*}\sideset{}{^*}\sum_{\substack{a\bmod p\\S(a,1;p)>0}}1&\geqslant\frac{1}{V_{2}(p)}\bigg(\sideset{}{^*}\sum_{\substack{a\bmod p\\S(a,1;p)>0}}S(a,1;p)\bigg)^2=\frac{4}{9\pi^2}p+O(p^{1/2}).\end{align*}

Similarly, we also have
\begin{align*}\sideset{}{^*}\sum_{\substack{a\bmod p\\S(a,1;p)<0}}1&\geqslant\frac{4}{9\pi^2}p+O(p^{1/2}).\end{align*}

\bigskip

\section{Additional remarks}

In a recent paper, D. I. Tolev \cite{To} obtained by elementary methods a smart identity of $S(a,b;p)$ for an arbitrary odd prime $p$ and $(ab,p)=1$. Namely,
\[S^2(a,b;p)=p+\sum_{x\bmod p}\left(\frac{x^2-4x}{p}\right)S(a,bx,p).\] From this identity, he can easily deduce that
\[|S(a,b;p)|\leqslant\sqrt{p+p^{3/2}},\ \ \ (ab,p)=1,\]which is of a smaller constant factor than the previous result of H. D. Kloosterman.

Of course, one may expect there exist certain identities for higher moments (especially for odd power moments) of Kloosterman sums. This still remains blank, however it seems a meaningful direction.

On the other hand, many arithmetical problems lead us to deal with the general Kloosterman sum twisted by a Dirichlet character $\chi\pmod p$, which is defined by
\[S_\chi(a,b;p)=\sideset{}{^*}\sum_{x\bmod p}\chi(x)e\left(\frac{ax+b\overline{x}}{p}\right).\] If $\chi$ is the principal character, this reduces to classical Kloosterman sum $S(a,b;p)$, if $\chi$ is the Legendre symbol mod $p$, this is known as the Sali\'{e} sum. Following Tolev's idea, one can find the following identity
\begin{align*}|S_\chi(a,b;p)|^2=p-1+\sideset{}{^*}\sum_{x\bmod p}S_\chi(x,x;p)e\left(-\frac{2x+ab\overline{x}}{p}\right)\end{align*}holds for any $\chi\bmod p$ and $(ab,p)=1.$ Notice that $S_\chi^2(a,b;p)=\chi(-\overline{a}b)|S_\chi(a,b;p)|^2,$ $(ab,p)=1$, thus a corresponding identity also holds for $S_\chi^2(a,b;p)$. Moreover, one can also deduce the following upper bounds that for $(ab,p)=1$,
\[|S_\chi(a,b;p)|<\begin{cases}
2^{1/4}p^{3/4},& \text{if }\chi=(\frac{\cdot}{p}),\\
p^{3/4},& \text{if }\chi\neq(\frac{\cdot}{p}).
\end{cases}\]

In closing, we would like to mention a twisted moment of Kloosterman sums studied by C. L. Liu \cite{L}. He proved that
\begin{align*}\sideset{}{^*}\sum_{a\bmod p}\chi(a)U_k(\cos\theta_p(a))\ll_kp^{1/2}\end{align*}
holds for any positive integer $k$ and odd prime $p$, where $\chi$ is a non-real character mod $p$. Clearly, his estimate yields
\[\sideset{}{^*}\sum_{a\bmod p}\chi(a)S^{2k}(a,1;p)\ll_k p^{k+1/2},\]which shows that there must be significant cancellations because of the existence of the character $\chi$. Regarding the case of odd powers, it follows that
\[\sideset{}{^*}\sum_{a\bmod p}\chi(a)S^{2k+1}(a,1;p)\ll_k p^{k+1}\]
and
\[\sideset{}{^*}\sum_{a\bmod p}\chi(a)|S^{2k+1}(a,1;p)|\ll_k p^{k+1}.\]
\bigskip

\noindent \textbf{Acknowledgements.} The present paper was prepared as YY participated in the 5th National Conference in Number Theory held in Zhaoqing. Both authors are grateful to Prof. C. L. Liu for his helpful suggestions during the conference, and to the referee for the valuable advice and comments. The work of the authors is supported by the Fundamental Research Funds for the Central Universities.
\bigskip

\bibliographystyle{plainnat}

\bigskip

\bigskip

\end{document}